\documentclass{amsart}

\usepackage{vmargin}

\newtheorem{theorem}{Theorem}[section] 
\newtheorem{lemma}[theorem]{Lemma}

\newtheorem{definition}[theorem]{Definition}

\newcommand{\Cay}{{\rm Cay}}

\newcommand{\ZZ}{\mathbb{Z}} 

\newcommand{\RR}{\mathbb{R}}

\newcommand{\diam}{\hbox{\rm
diam}} \newcommand{\Diam}{\hbox{\rm Diam}}

\newcommand{\Label}{\label}
 
\newcommand{\la}{\langle}
\newcommand{\ra}{\rangle} 
\renewcommand{\P}{{\mathcal{P}}}

\newcommand{\B}{{\mathcal{B}}}

\newcommand{\dist}{\partial} \newcommand{\co}[1]{\overline{#1}}

\newcommand{\gf}{G}

\title{The diameter of random Cayley digraphs of given degree} \date{\today}

\author{Manuel E. Lladser} \address{Department of Applied Mathematics,
University of Colorado, Boulder, CO 80309-0526, THE UNITED STATES}
\email{manuel.lladser@colorado.edu}

\author{Primo\v{z} Poto\v{c}nik} \address{Faculty of Mathematics and
Physics, University of Ljubljana, Ljubljana, SLOVENIA}
\email{primoz.potocnik@fmf.uni-lj.si}

\author{Jozef \v{S}ir\'a\v{n}} \address{Department of Mathematics,
University of Auckland, Private Bag 92019, NEW ZEALAND}
\email{siran@math.auckland.ac.nz}

\author{Jana \v{S}iagiov\'a} \address{Department of Mathematics, SvF,
Slovak University of Technology, Bratislava, SLOVAKIA} \email{}

\author{Mark C. Wilson} \address{Department of Computer Science, University
of Auckland, Private Bag 92019 Auckland, NEW ZEALAND}
\email{mcw@cs.auckland.ac.nz}

\begin{document} 

\begin{abstract}

We consider random Cayley digraphs of order $n$ with uniformly
distributed generating set of size $k$. Specifically, we are interested
in the asymptotics of the probability such a Cayley digraph has diameter
two as $n\to\infty$ and $k=f(n)$. We find a sharp phase transition from
0 to 1  as the order of growth of $f(n)$ increases past $\sqrt{n \log
n}$. In particular, if $f(n)$ is asymptotically linear in $n$, the
probability converges exponentially fast to $1$.
\end{abstract}

\maketitle

\section{Introduction} \Label{sec:intro}

It is well known that almost all graphs and digraphs have diameter two
\cite{bollobas-graph-theory}. This result has been generalized and
strengthened in various directions, of which we shall be interested in
restrictions to Cayley graphs and digraphs.

In \cite{meng-liu} it was proved that almost all Cayley digraphs
have diameter two, and in \cite{meng-huang} this was extended to
Cayley graphs. The random model used in \cite{meng-liu,
meng-huang} is the most straightforward one: in terms of Cayley
digraphs for a given group $G$, one chooses a random generating
set by choosing its elements among the non-identity elements of
$G$ independently and uniformly, each with probability $2^{-n+1}$
where $n$ is the order of $G$. Observe that such generating sets
have size at least $n/2$ with probability at least $1/2$, in which
case the corresponding Cayley digraphs automatically have diameter
at most two. The less trivial part of \cite{meng-liu} therefore
concerns random Cayley digraphs in which the number of generators
is at most half of the order of the group.

This motivates a study of random Cayley digraphs in which the
number of generators is restricted. The fundamental problem here is
the following: for which functions $f$ is it true that the diameter
of a random Cayley digraph of an arbitrary group of order $n$ and
of degree $f(n)$ is asymptotically almost surely equal to $2$ as
$n$ tends to infinity? By the well known Moore bound for graphs or
digraphs of diameter two we know that $f$ has to increase at least
as fast as $\sqrt{n}$. However, even the case when $f(n)=cn$ for a
constant $c$ seems not to have been investigated before and, as we
shall see, leads to interesting questions in the study of
generating functions.

In order to investigate the above problem one cannot use the model of
\cite{meng-liu}. Instead, we will consider the uniform distribution of
subsets of size $k$ in the set of all non-identity elements of a given
group of order $n$. A detailed description of the model and the
associated parameters is given in Section~\ref{sec:model}. The
probability that a random Cayley digraph of (in- and out-) degree $k$ on
a group of order $n$ has diameter $2$ will be estimated in
Section~\ref{sec:estimates} in terms of a certain combinatorial function
$p(n,k,t)$ where $t$ is a parameter that depends on the group and
$2t<n$. Dependence on the group is then eliminated by showing that one
can use $t=\lfloor \gamma n \rfloor$ for a suitable constant $\gamma \le
1/2$ in the estimates. Setting $k=f(n)$ and $t= \lfloor\gamma n\rfloor$,
the probability that a random Cayley digraph has diameter two can be
studied by means of the asymptotic behaviour of $p(n,f(n),\lfloor\gamma
n\rfloor)$ as $n\to \infty$. The two cases of particular interest are
$f(n)=\lfloor cn \rfloor $ for a fixed constant $c$ with $0< c <1/2$, and
$f(n)=\lfloor n^{\alpha} \rfloor$ for a fixed constant $\alpha$ such that $1/2 < \alpha
< 1$. By a delicate asymptotic analysis, in Section~\ref{sec:GF} we
prove that in both cases the diameter of a random Cayley graph is
asymptotically almost surely equal to two. We also consider analogous
questions for random Cayley graphs on elementary abelian $2$-groups.
Under this restriction we obtain tighter bounds in terms of $p(n, k, t)$ for the
probabilities, which raises an interesting question on the probability
evolution if $f(n)\sim c\sqrt{n}$ for a constant $c>1$.

\section{The model} \Label{sec:model}

Throughout, let $G$ be a finite group of order $n$ and let $k$ be a
positive integer not exceeding $n-1$. The set of non-trivial elements of
$G$ will be denoted by $G^*$. For a set $A$ and an integer $r$, the symbol
$\binom{A}{r}$ will stand for the set of all subsets of $A$ of size $r$.

For  $S\in {G^*\choose k}$, the {\em Cayley digraph on $G$ relative to
$S$}, denoted by $\Cay(G,S)$, is the $k$-valent digraph with vertex set $G$
and arc set $\{(g,gs) : g \in G, s \in S\}$. The {\em distance}
$\dist(g,h)$ from the vertex $g$ to the vertex $h$ in $\Cay(G,S)$ is the
length of the shortest directed path from $g$ to $h$ in $\Cay(G,S)$. The
{\em diameter} $\diam(\Cay(G,S))$ is the smallest integer $d$ such that for
every ordered pair $(g,h)$ the distance from $g$ to $h$ is at most $d$.

We are now ready to introduce our model for random Cayley digraphs of a
given valence. Let $\P(G,k)$ be the probability space $(\B,2^{\B},P)$ where
$\B={G^* \choose k}$, $2^{\B}$ is the power set of $\B$, and $P$ is the
uniformly distributed probability measure on $\B$. Since $|\B| =
{n-1\choose k}$, a simple counting argument shows that $\Pr(\{S\}) =
{n-1\choose k}^{-1}$ for every $S\in \B$. More generally, for every subset
$L\subseteq G^*$ of size $\ell$, the probability that a random set $S\in\B$
contains $L$ as a subset is given by

\begin{equation} \label{eqn:1} \Pr( S\supseteq L) =  \Pr(\{ S \in {G^*\choose
k} : L\subseteq S\}) = {{n-1-\ell}\choose{k-\ell}}{n-1\choose
k}^{-1}=\frac{(k)_{\ell}}{(n-1)_{\ell}} \end{equation}

\noindent where $(r)_{\ell}=r(r-1)\ldots (r-\ell+1)$ denotes the $\ell$-th
descending factorial of $r$ (with the convention that $r_0=1$). We can now
define a random variable $\Diam \colon \B \to \RR$ on the probability space
$\P(G,k)$ by letting, for every $S \in {G^*\choose k}$,

\begin{equation}
\label{eqn:2}
\Diam(S) = \diam(\Cay(G,S)).
\end{equation}

\noindent The main goal of this article is to derive bounds on the
probability of the event $\{S \in {G^*\choose k} :
\diam(\Cay(G,S)) = 2\}$ and study the asymptotic behaviour of the bounds.

Since Cayley digraphs are vertex-transitive, the diameter of $\Cay(G,S)$
coincides with the maximum value of $\dist(1,y)$ over all $y\in G^*$.
Clearly, if $\dist(1,y) \le 2$, then $y\in S$, or there exists $x\in S$
such that $(1,x,y)$ is a directed path from $1$ to $y$ of length $2$. The
latter is equivalent to requiring that $\{x,x^{-1}y\} \subseteq S$. This
shows that the following events will play an important role:

\begin{definition} For $x,y\in G^*$, let
$$
T(x,y) = \{ S : S \in \binom{G^*}{k},\, \{x,x^{-1}y\} \subseteq S\}
\quad \hbox{ and } \quad X(y) = \bigcup_{x\in G^*} T(x,y).
$$
\end{definition}

Let $S$ be an arbitrary element of ${ G^* \choose k}$. Clearly, there is a
directed path  from $1$ to $y$ of length $2$ in $\Cay(G,S)$ if and only if
$S\in X(y)$. In other words, $S\in \co{X(y)}$ if and only if there is no
directed path from $1$ to $y$ in $\Cay(G,S)$ of length exactly $2$. Thus,
if $\diam (\Cay(G,S))> 2$ then $S\in
\cup_{y\in G^*} \co{X(y)}$. Therefore we have the following inequality:

\begin{equation}
\label{eqn:3} \Pr(\Diam > 2) \> \le \> \sum_{y\in G^*} \Pr(\co{X(y)}).
\end{equation}

On the other hand, if $\diam(\Cay(G,S)) \le 2$, then for every $y\in
G^*$ we have $y\in S$ or $S\in X(y)$. Hence $\Pr(\Diam \le 2) \le
\Pr(X(y)) + \Pr(y \in S)$, and by \eqref{eqn:1}, it follows that
$\Pr(\Diam \le 2) \le \Pr(X(y)) + \frac{k}{n-1}$, which is equivalent to
$\Pr(\Diam > 2) \ge \Pr(\co{X(y)}) -
\frac{k}{n-1}$. This, together with \eqref{eqn:3}, shows that

\begin{equation} \label{eqn:4} 
M - \frac{k}{n-1} \> \le \> \Pr(\Diam > 2) \>
\le \> (n-1) M, \quad \hbox{ where }\> M = \max_{y\in G^*} \Pr(\co{X(y)}).
\end{equation}

The inequality \eqref{eqn:4} provides the basis for our
investigation. In what follows we consider estimates for the
quantity $M$ appearing in \eqref{eqn:4}.

\section{The estimates}
\label{sec:estimates}

The key to deriving bounds on $M$ is the evaluation of the probability
$\Pr(\co{X(y)}) = 1- \Pr(\cup_{x\in G^*} T(x,y))$. As Lemma~\ref{lem:1} below
shows, this probability is closely related to the values of 
$p(n,k,t)$, $t,k\le n$, where $p(n, k, t)$ is defined by

\begin{equation} 
\label{eqn:F} p(n,k,t) = \sum_{i= 0}^t(-1)^{i}{t \choose
i}{{{n-1-2i}\choose {k-2i}}}{n-1\choose k}^{-1} =\sum_{i= 0}^t(-1)^{i}{t
\choose i}\frac{(k)_{2i}}{(n-1)_{2i}}. 
\end{equation}

\begin{lemma} \label{lem:1} Let $y\in G^*$ and let $J\subseteq
G^*\setminus \{y\}$ be a set of size $t$ such that the sets
$\{x,x^{-1}y\}$, with $x\in J$, are pairwise disjoint and of size $2$.
Then $$ \Pr(\co{X(y)}) \le  p(n,k,t). $$ 
\end{lemma}

\noindent \begin{proof} We start with a simple inequality
$$\Pr(\co{X(y)}) = 1 - \Pr(X(y)) = 1 - \Pr(\cup_{x\in G^*}T(x,y)) \le 1
- \Pr(\cup_{x\in J}T(x,y)).$$ Now observe that the set $\cap_{x\in
I}T(x,y)$ consists of all those $S\in {G^*\choose k}$ for which
$\cup_{x\in I}\{x,x^{-1}y\} \subseteq S$. Hence, if $I\subseteq J$
and $|I|=i$, then 
$$ \Pr(\cap_{x\in I}T(x,y))=
{{n-1-2i}\choose{k-2i}}{n-1\choose
k}^{-1}=\frac{(k)_{2i}}{(n-1)_{2i}}. $$ 
By the inclusion-exclusion
formula, we have  
\begin{align*} \Pr(\cup_{x\in J}T(x,y)) &=
\sum_{i=1}^{t} (-1)^{i-1}\sum_{I\subseteq {J\choose
i}}\Pr(\cap_{x\in I}T(x,y)) \\ & =\sum_{i= 1}^t(-1)^{i-1}{t \choose
i}{{{n-1-2i}\choose {k-2i}}}{n-1\choose k}^{-1}\\ & = \sum_{i=
1}^t(-1)^{i-1}{t \choose i}\frac{(k)_{2i}}{(n-1)_{2i}};
\end{align*} and the result follows. 
\end{proof}

\medskip A straightforward consequence of the above proof is the fact that
$0\le p(n,k,t)\le 1$ and that $p(n, k, t)$ is decreasing in $t$ (in the range for
which there is an appropriate group $G$ for which the above lemma can be
used).

We continue with establishing an upper bound on the parameter $t$
appearing in the sums above. For this we need an estimate of the
number of `square roots' of a non-identity element in a group.
Therefore, for an element $y\in G$ let $\sigma(y)$ denote the set
of $x\in G$ such that $x^2 =y$.

\begin{lemma}
\label{lem:sigma}
If $y$ is a non-trivial element of a finite group $G$, then
$|\sigma(y)| \le \frac{3}{4}|G|$.
\end{lemma}

\begin{proof}
Suppose the contrary and let $y \in G^*$ be such that $|\sigma(y)|
> \frac{3}{4}|G|$. Take an arbitrary element $z\in G$, and observe
that $\sigma(y^z)=\sigma(y)^z$. In particular, $|\sigma(y^z)| =
|\sigma(y)| > \frac{3}{4}|G|$, implying that $\sigma(y) \cap
\sigma(y^z)$ is non-empty. Take any $x\in\sigma(y) \cap
\sigma(y^z)$, and note that $y=x^2=y^z$. Hence $y$ is in the
centre of $G$. Now consider the quotient projection $\pi \colon G
\to G/\la y\ra$, let $s$ be the order of $y$ in $G$, and let $T$
denote the set of elements $x\in G/\la y\ra$ such that $x^2 =1$.
Clearly, $\pi(\sigma(y)) \subseteq T$. Suppose that
$\pi(x_1)=\pi(x_2)$ for some pair of $x_1, x_2 \in \sigma(y)$,
$x_1,\not = x_2$. Then $x_2=x_1y^r$ for $1\le r < s$, and so
$y=x_2^2 = x_1^2 y^{2r} = y^{1+2r}$, implying that $s=2r$. This
shows that the $\pi$-preimage of each element in $T$ contains at
most $2$ elements from $\sigma(y)$, and contains at most one
element from $\sigma(y)$ if $s$ is odd. Therefore $|T| \ge
\frac{|\sigma(y)|}{2} > \frac{3s}{8}|G/\la y \ra|$ if $s$ is even,
and $|T| \ge |\sigma(y)| > \frac{3s}{4}|G/\la y \ra|$ if $s$ is
odd. On the other hand, $|T| \le |G/\la y \ra|$, implying that
$s=2$ and $|T| > \frac{3}{4}|G/\la y \ra|$. It is known that the
only groups for which the proportion of the involutions is more
than $\frac{3}{4}$ are the elementary abelian $2$-groups (see
\cite{wall-involutions}). Hence $G/\la y \ra$ is an elementary
abelian $2$ group, say $G/\la y \ra \cong \ZZ_2^p$, where $p$ is a
positive integer, and consequently, $G\cong \ZZ_2^{p+1}$ or
$G\cong \ZZ_2^{p-1} \times \ZZ_4$.  However, it is clear that no
element $y$ in such groups satisfies $|\sigma(y)| >
\frac{3}{4}|G|$.
\end{proof}

\medskip We remark that the bound in the previous lemma is sharp. For
example, if $Q=\{\pm 1, \pm i, \pm j, \pm k\}$ is the quaternion group and
if $G=Q\times \ZZ_2^q$, then for the element $y=(-1, 0, \ldots, 0)$ we have
$\sigma(y)= \frac34 |G|$.

We are now ready to prove the main result of this section, which is an
upper bound on the probability $\Pr(\Diam > 2)$ in terms of $p(n,k,t)$.

\begin{theorem} \label{the:1} Let $G$ be a finite group, and let $k$ be
such that $1\le k\le n=|G|$. Then for the random variable {\rm Diam} on the
probability space $\P(G,k)$ we have

$$
\Pr(\Diam > 2) \le (n-1)p(n,k,\lfloor (n-4)/12 \rfloor).
$$
\end{theorem}

\begin{proof}
Let $y\in G^*$ and let  $s=|\sigma(y)|$. We first show that there
exists a set $J\subseteq G\setminus \{1,y\}$ of size at least
$t=\lfloor\frac{n-1-s}{3}\rfloor$ such that the sets $\{x,x^{-1}y\}$, where
$x\in J$, are pairwise disjoint and of size $2$. We shall define such a set
$J$ recursively.

If $s \ge n-3$, then $t = 0$, and $J=\emptyset$ will do the job. So we may
assume that $s \le n-4$. Then the set $C=G^*\setminus (\sigma(y) \cup \{y
\})$ is non-empty, and we can choose $x_1\in C$ and set $J_1=\{x_1\}$.

Now suppose that $J_1, \ldots, J_{\ell}$ have been already defined for some
$\ell < t$, and suppose that this has been done in such a way that
$J_i\subseteq C$, $|J_i|=i$ and $|\cup_{x\in J_i}\{x,x^{-1}y\}| = 2i$ for
every $i\in \{1,\ldots, \ell\}$. Let $K_\ell = \cup_{x\in J_{\ell}} \{x,
yx^{-1},x^{-1}y\}$. Then $|K_\ell|\le 3\ell \le 3t -3 \le n-s -4 < n-s -2=
|C|$. Hence, we can choose an element $x_{\ell+1}\in C \setminus K_\ell$,
and define $J_{\ell +1} = J_\ell \cup \{x_{\ell+1}\}$. Clearly,
$|J_{\ell+1}| = |J_\ell|+1 = \ell +1$, and since $x_{\ell+1}\in C$, also
$J_{\ell+1} \subseteq C$. Now suppose that $|\cup_{x\in
J_{\ell+1}}\{x,x^{-1}y\}| < 2\ell + 2$. Then one of the elements
$x_{\ell+1}$ and $x_{\ell+1}^{-1}y$ belongs to $\cup_{x\in
J_\ell}\{x,x^{-1}y\}$. However, both cases imply that $x_{\ell+1} \in
K_\ell$, which contradicts our assumption. Hence this construction yields a
set $J=J_t$ with the desired properties.

From Lemma ~\ref{lem:1} it follows that $\Pr(\co{X(y)}) \le p(n,k,\lfloor
(n-1-s)/3\rfloor)$. Then by Lemma \ref{lem:sigma} we have $s\le \frac34 n$
and therefore $(n-1-s)/3\ge (n-4)/12$. Using \eqref{eqn:4}, and the fact
that $p(n, k, t)$ is decreasing  in $t$ we arrive at the
inequality in the statement of the theorem.
\end{proof}

\medskip If the group $G$ is an elementary abelian $2$-group, then the
value of $M$ in \eqref{eqn:4} can be expressed in terms of $p(n,k,t)$ exactly.
This leads to the following result.

\begin{theorem} \label{the:2} Let $G\cong \ZZ_2^d$ be an elementary abelian
$2$-group, and let $1\le k\le n=2^d$. Then for the random variable {\rm
Diam} on the probability space $\P(G,k)$ we have

$$
p(n,k,(n-2)/2)-\frac{k}{n-1}\le \Pr(\Diam > 2) \le (n-1)p(n,k,(n-2)/2 ).
$$
\end{theorem}

\begin{proof} We show that for any $y\in G^*$ we have $\Pr(\co{X(y)}) =
p(n,k,\frac{n-2}{2})$. Let $J_y$ be any transversal of the
subgroup $\la y \ra$ in $G$, and let $J = J_y \setminus \{1,y\}$.
Then $\Pr(X(y)) = \Pr(\cup_{x\in G^*}T(x,y)) = \Pr(\cup_{x\in
J}T(x,y))$. On the other hand, the set $J$ satisfies the
conditions of Lemma~\ref{lem:1}. Since $|J| = \frac{n-2}{2}$, we
have $$ \Pr(\co{X(y)}) = 1- \Pr(\cup_{x\in J}T(x,y)) = p(n,k,(n-2)/2).
$$ The statement now follows from \eqref{eqn:4}. \end{proof}

It is now clear that knowledge of the asymptotic behaviour of $p(n,k,t)$
would allow us to make conclusions about the asymptotic behaviour
of the random variable $\Diam$. As explained in the Introduction,
the most interesting cases to study are $f(n)=\lfloor cn\rfloor$ for $0<c<1/2$
and $f(n)=\lfloor n^{\alpha}\rfloor$ for $1/2<\alpha<1$. For example, Theorem
\ref{the:1} shows that if $\lim_{n\to \infty}(n-1)p(n,\lfloor
cn\rfloor,\lfloor (n-4)/12\rfloor)=0$ for a constant $c$ such
that $0<c<1/2$, then the diameter of a random Cayley digraph of
order $n$ and degree $\lfloor cn\rfloor$ is asymptotically almost
surely equal to two. By the same token, if $\lim_{n\to
\infty}(n-1)p(n,\lfloor n^{\alpha}\rfloor,\lfloor (n-4)/12
\rfloor)=0$ for $1/2<\alpha < 1$, then the diameter of a random
Cayley digraph of order $n$ and degree $\lfloor n^{\alpha}\rfloor$
is also asymptotically almost surely equal to two. Similar 
statements, with $n$ a power of two and with $(n-2)/2$ in the
third position, hold for random abelian Cayley digraphs
on the basis of Theorem \ref{the:2}. In the next section we will
show that the above limits are indeed equal to zero and therefore
the corresponding random Cayley digraphs almost surely have
diameter two. Since Theorem \ref{the:2} gives also a lower bound,
it is natural to ask if, for $n$ a power of two, $\lim_{n\to
\infty}p(n,\lfloor cn^{1/2}\rfloor, (n-2)/2)=1$ for sufficiently large 
$c$. As we shall see, the answer to this question is in the
negative. In what follows we also describe more precisely the
threshold at which $\Pr(\Diam \leq 2)$ jumps asymptotically away
from $0$.

\section{Asymptotic analysis} \label{sec:GF}

We use a mixture of techniques, based on generating functions, with
varying levels of sophistication. Some of the questions of the previous
section are quickly addressed by relatively simple means, while for
others it is cleaner to apply asymptotic techniques for the analysis of
coefficients of multivariate generating functions as developed in
\cite{p1, p2, p5, baryshnikov-pemantle, lladser03}. See \cite{p9} for a
detailed survey of the use of such techniques in combinatorial problems.
For the hardest questions we use the recently developed machinery of
\cite{lladser07a, lladser06}.

The quantity $a(n,k,t) := \binom{n}{k} p(n+1,k,t)$ is simpler to analyse
in this way than $p(n,k,t)$ itself. It is easily seen from above that
$a(n, k, t)$ has a purely combinatorial description. Namely, given a set
of size $n$, we choose $t$ disjoint pairs from this set. Then $a(n,k,t)$
is the number of subsets of size $k$ that contain none of the pairs.
Note  that $a(n, k, t) = 0$ if $k + t > n$, by the pigeonhole
principle (since the complement of $S$ has size less than $t$, $S$ must
contain at least $t+1$ of the $2t$ paired elements).

From the statement of the problem, $a(n, k, t)$ is not defined if $2t >
n$; however, formula~\eqref{eqn:F} still makes sense in that case, even
though it does not define the probability of any event. In fact, $a(n,
k, t)$ can be negative for large $t$. Asymptotic analysis in this case
is considerably more difficult than what is presented below, and we will
avoid this case in the present paper, since it is not relevant to the
original combinatorial question.

\subsection{Generating functions}
\label{ss:gf}

We first compute the trivariate generating function of $a(n,k,t)$. The
most direct approach is to use some well-known bivariate generating
functions $\sum_{i,j} a_{ij} x^i y^j$. Throughout, we use the convention
that the binomial coefficient $\binom{k}{l}$, with $k, l \in
\mathbb{Z}$, is zero unless $0 \leq l \leq k$. If $a_{ij} =
\binom{i+j}{i}$ then the generating function is $(1- x - y)^{-1}$, while
that for $a_{ij} = \binom{i}{j}$ is $(1- x(1+y))^{-1}$. We now compute
\begin{align*} 
\sum_{n, k, t, i}  x^n y^k z^t w^i \binom{t}{i}
\binom{n - 2i}{k - 2i} & = \sum_{N, K, i, j} x^{N+2i} y^{K+2i}
z^{i+j} w^i \binom{i+j}{i} \binom{N}{K} \\ & =  \left(\sum_{i, j}
\binom{i+j}{i}  (zwx^2y^2)^i z^j \right) \left(\sum_{N, K}
\binom{N}{K} x^N y^K \right) \\ & = \frac{1}{1 - z(1 + wx^2y^2)}
\frac{1}{1 - x(1 + y)} ,
\end{align*} 
which yields the trivariate generating function
\begin{equation} 
\gf(x,y,z) = \sum_{n,k,t} a(n,k,t)
x^n y^k z^t = \frac{1}{1 - z(1 - x^2 y^2)} \frac{1}{1 - x(1+y)} \, .
\end{equation}

Note that if we impose the restriction $2t \leq n$, then the sum over $N$ is restricted to $N \geq 2j$. Now summing over $N, K, i, j$ as above we
obtain the more relevant restricted trivariate generating function
\begin{equation} 
\gf_1(x, y, z) = \sum_{\{n,k,t \colon 2t
\leq n\}} a(n,k,t) x^n y^k z^t =  \frac{1}{1 - x(1+y)}\frac{1}{1 -
zx^2(1+2y)} =: \frac{1}{H_1 H_2}. 
\end{equation}

The series $\gf_1$ is more useful for our purposes, since all coefficients
are nonnegative.

\subsection{Basic asymptotic approximations}
\label{sec:asymp-basics}

We list here some standard asymptotic approximations that will be used later.

\begin{lemma}
Write $n = \lambda k$ with $0 < \lambda < 1$. Then
\begin{equation}
\label{eq:binom-approx-linear}
\binom{n}{k} = \exp (n R(\lambda)) \, P(\lambda) \, n^{-1/2} \, C(n, \lambda)
\end{equation}
where
\begin{align*}
R(\lambda) & = - \lambda \log \lambda - 
(1 - \lambda) \log (1 - \lambda) \\
P(\lambda) & = \left(2\pi \lambda (1 - \lambda)\right)^{-1/2} \\
C(n, \lambda) & = (1 + O(n^{-1}) + O((n\lambda)^{-1}) + O((n(1-\lambda))^{-1}) 
\qquad \text{as $n \to \infty$}.
\end{align*}

\end{lemma}

\begin{proof}
A direct application of Stirling's approximation.
\end{proof}

We call $R$ the exponential rate and $P$ the leading coefficient, while
$C$ is the correction term. Of course $R$ depends on $\lambda$ and so
may vary with $n$ as $n \to \infty$.

\begin{lemma}
For $t \geq k \geq 0$ define
$$
b(t,k):=\frac{2^k \binom{t}{k}}{\binom{2t}{k}}.
$$
Then with $t = \lambda k$ we have
\begin{equation}
\label{eq:b-approx-linear}
b(t, k) = \exp(t R(\lambda))\,  P(\lambda)\, C(t, k)
\end{equation}
where
\begin{align*}
R(\lambda) &:= (2 - \lambda) \log (1 - \lambda/2) - (1 - \lambda) \log (1 - \lambda)) \\
P(\lambda) &:= \left(\frac{2 - \lambda}{2 - 2\lambda}\right)^{1/2} \\
C(\lambda) &:= 1 + O(t^{-1}) + O((t \lambda)^{-1}) + O((t(1 - \lambda))^{-1}) \qquad \text{as $t \to \infty$}.
\end{align*}
\end{lemma}

\begin{proof}
An immediate application of the previous lemma (replace $n$ by $2t$ and $\lambda$ by 
$\lambda/2$ in the denominator and replace $n$ by $t$ in the numerator).
\end{proof}

We also need some standard facts about the stationary phase
approximation of an oscillatory integral. We recall them below and refer
to a standard text such as \cite{wong} for details. Define
$$
I(f; n) = \int_{a}^b e^{n f(\theta)} g(\theta) \, d\theta
$$
where $f$ and $g$ are smooth functions and $\Re f \geq 0$ on $[a, b]$.

\begin{lemma}[Laplace approximation]
\label{thm:Laplace}

Let $I(f; n)$ be as defined above. Suppose that $\Re f > 0$ on $[a, b]$
except at a single point $x \in (a, b)$. Furthermore suppose that $f'(x) = 0$ and
$f''(x) \neq 0$, while $g(x) \neq 0$. Then 
$$
I(f;n) = \exp(nf(x)) \frac{g(x)}{\sqrt{2\pi n f''(x)}} (1 + O(n^{-1})) \qquad 
\text{as $n \to \infty$}.
$$
The implied constant in the $O$-term remains bounded as we vary $f$ and
$g$ provided that no hypotheses change, $x$ remains in a compact subset
of $(a, b)$, $f''(x)$ remains bounded away from zero, and the maximum of
$|g|$ remains bounded.

\hfill $\Box$
\end{lemma}

\subsection{Abelian groups}
\label{ss:abelian}

\subsubsection{The linear case}
\label{sss:abelian-linear}

The simpler form of the bounds involving $p(n,k,t)$ in the abelian case allows an easy 
elementary approach, which we now present.

Suppose that $t = (n-2)/2$, so that $n - 1 = 2t + 1$. A direct evaluation as 
in Section~\ref{ss:gf} above shows that 
$$ 
\sum_{k,t} a(2t+1,k,t) y^k z^t  = \frac{1+y}{1 - z(1 + 2y)} 
$$
and hence we may extract coefficients to obtain
$$ 
a(2t+1, k, t) = 2^k \binom{t}{k} + 2^{k-1} \binom{t}{k-1}. 
$$
Now 
\begin{align*} 
p(2t+2, k, t) & = \binom{2t+1}{k}^{-1} a(2t+1, k, t) \\
& = \left[\binom{2t}{k}^{-1} 2^k \binom{t}{k} \right] \left[ \frac{(2t - k
+ 2)(2t + 1 - k)}{(2t - 2k + 2) (2t + 1)} \right] \\ 
&=: [b(t,k)][c(t,k)].
\end{align*} 
When $ k = t + 1$, the right side above should be replaced  by $2^t$, but we do not deal
with this case below anyway (since $k < n/2$, $n$ is even and $k$ is an
integer we must have $k \leq (n - 2)/2 = t)$.

From \eqref{eq:b-approx-linear}, we obtain the exponential rate of $b(t,
k)$ with respect to $t$ as 
$$
R(\lambda) = (2 - \lambda) \log(1 - \lambda/2) - (1 - \lambda)\log(1 - \lambda)
$$
with $\lambda = k/t$. This is easily
seen by elementary calculus to be negative for $0 < \lambda \leq 1$.
Furthermore $c(t, k)$ has exponential rate zero. Thus in combination
with Theorem \ref{the:2}  we have:

\begin{theorem}\label{the:4}
For any constant $c$ such that $0 < c < 1/2$, the diameter of a random
Cayley digraph on an elementary abelian $2$-group of order $n$ and
degree $\lfloor cn\rfloor$ is asymptotically almost surely equal
to two. Furthermore the convergence is exponentially fast.
\hfill $\Box$
\end{theorem}

\subsubsection{The sublinear case}
\label{sss:abelian-sublinear}

We now consider the case where $k$ is of order $n^{\alpha}$ with $1/2 <
\alpha < 1$. For $k = \lambda t$ with $\lambda = o(1)$ as $t \to
\infty$, we have $c(t, k) = 1 + O(\lambda)$. By
\eqref{eq:b-approx-linear} we again have 
$$
R(\lambda) =  (2 - \lambda) \log(1 - \lambda/2) - (1 - \lambda)
\log(1 - \lambda) = -\lambda^2/4 + O(\lambda^3).
$$
Thus if $k$ grows at least as fast as $n^{\alpha}$ with $\alpha > 1/2$,
it is definitely the case that the upper bound $(2t+1) b(t, k) c(t, k)$
decays faster than polynomially. Using Theorem \ref{the:2} again, we
have the following conclusion.

\begin{theorem}\label{the:5}
For any constant $\alpha$ such that $1/2<\alpha<1$, the diameter
of a random Cayley digraph on an elementary abelian group of order
$n$ and degree $\lfloor n^{\alpha}\rfloor$ is asymptotically
almost surely equal to two. 
\hfill $\Box$
\end{theorem}

Also, the approximation above shows that if $k = \lfloor c
\sqrt{n}\rfloor$ then the lower bound $b(t, k) c(t, k)$ converges to
$\exp(-c^2/2)$, and not $1$. This shows that for $n$ a power of two,
$p(n,\lfloor cn^{1/2}\rfloor, (n-2)/2)$ does {\em not} tend to $1$ as
$n\to\infty$. 

The case of general groups could be attacked in a similar way to the
elementary approach above, but with more effort. For example, $p(n,k,t)$ is
decreasing in $t$, so that $p(n, k, (n - 2)/12) \geq p(n, k, \lfloor
(n-4)/12 \rfloor) \geq p(n, k, (n-4)/12)$ for sufficiently large $n$. As
above we can compute the bivariate generating function for $F(12t + 3,
k, t)$ and $F(12t + 1, k, t)$, and estimate each of these as above when
$k$ is of order $n^{\alpha}$. However, this approach depends heavily on
the relatively nice formula for $p(n,k,t)$ involving well-studied binomial
coefficients. For variety, and to illustrate that more detailed
expansions can be obtained in more generality, we use a different
approach in Section~\ref{sec:parameter-varying}.

\subsection{General groups}
\label{sec:parameter-varying}

In this section we use Theorem~\ref{the:1} to study the asymptotic
behavior of the diameter of a general random Cayley digraph of order $n$
and degree $k$. We again consider two different regimes, namely $k=\lfloor c
n\rfloor$ and $k=\lfloor n^{\alpha}\rfloor$ where $0<c<1/2<\alpha<1$. 

Our analysis is in terms of parameter-varying integrals
that lead to uniform asymptotic expansions for the coefficients
$a(n,k,t)$, with $t=\lfloor(n-4)/12\rfloor$. We make heavy use of the
fact that $n\ge 2k$ and $n\ge2t$ for both regimes to reduce the problem
to the analysis of a one-dimensional parameter-varying integral. For the
first regime, the asymptotic behavior of the resulting integral relies
on the stationary phase method as found for example in~\cite{wong}.
Instead, for the second regime, the analysis follows the lines
of~\cite{lladser07a} and~\cite{lladser06} to properly use the Laplace approximation of Lemma~\ref{thm:Laplace}.

In what follows it is always assumed that $n\ge2k\ge0$ and $n\ge2t\ge0$.
We first reduce the problem to computing the asymptotics of a certain one-dimensional 
integral. Since $a(n,k,t)=[x^ny^kz^t]\gf_1(x,y,z)$, we obtain 
\begin{align*}
a(n,k,t)&=[x^ny^kz^t]\sum_{l=0}^\infty\frac{z^lx^{2l}(1+2y)^l}{1-x(1+y)},\\
&=[x^{n-2t}y^k]\frac{(1+2y)^t}{1-x(1+y)},\\
&=[x^{n-2t}y^k]\sum_{l=0}^\infty x^l(1+y)^l(1+2y)^t,\\
&=[y^k](1+y)^{n-2t}(1+2y)^t.
\end{align*}
Using Cauchy's formula, the above implies for all $r>0$ that
\begin{align}
\label{ide:coeff G2}
a(n,k,t)&=\frac{r^{-k}}{2\pi}\int_{-\pi}^\pi(1+re^{i\theta})^{n-2t}(1+2re^{i\theta})^t
e^{-ik\theta}d\theta.
\end{align}
In particular, we can rewrite
\begin{align}
\label{eq:fact a}
a(n,k,t)=(2\pi)^{-1} \cdot E(r;n,k,t)\cdot I(r;n,k,t),
\end{align}
where
\begin{align*}
E(r;n,k,t)&:= r^{-k}(1+r)^{n-2t}(1+2r)^t,\\
I(r;n,k,t)&:=\int_{-\pi}^\pi\left(\frac{1+re^{i\theta}}{1+r}\right)^{n-2t}\left(\frac{1+2re^{i\theta}}{1+2r}\right)^te^{-
ik\theta}d\theta.
\end{align*}

The integral in \eqref{ide:coeff G2} has been normalized by the factor
$(1+r)^{n-2t}(1+2r)^t$ to emphasize that the modulus of each of the two
factors in the integrand is maximized at $\theta=0$. 

To determine the asymptotic behavior of $a(n,k,t)$ the goal is to tune
$r$ with $(n,k,t)$ so that $I(r;n,k,t)$ decays polynomially with $n$, in other words
so that $E(r;n,k,t)$ captures the precise exponential growth rate
of the coefficients $a(n,k,t)$.

To accomplish our goal, motivated by the
stationary phase method, we rewrite the integrand of $I(r;n,k,t)$ in an
exponential-logarithmic form to obtain

\begin{align}
\label{eq:def I(r;n,k,t)}
I(r;n,k,t)&=\int_{-\pi}^\pi \exp\big\{-n\cdot F(\theta;r,d_1, d_2, d_3)\big\}d\theta,
\end{align}
where
\begin{eqnarray*}
\label{eq:def F}
F(\theta;r,d_1,d_2,d_3)&:=&
d_3\cdot i\theta-d_1\cdot\ln\left\{\frac{1+re^{i\theta}}{1+r}\right\}-d_2
\cdot\ln\left\{\frac{1+2re^{i\theta}}{1+2r}\right\},\\
d_1&:=& \frac{n-2t}{n},\\
d_2&:=&\frac{t}{n},\\
d_3&:=&\frac{k}{n}.
\end{eqnarray*}
In what follows all logarithms are to be interpreted in the principal
sense. In addition, unless otherwise stated, $d_1$, $d_2$ and $d_3$
always stand as short forms of the functions defined above. As a note on
our terminology, we refer to $I(r;n,k,t)$ as a parameter-varying
integral because $F(\theta;r,d_1, d_2, d_3)$, the so called phase term,
depends upon the parameter $n$ itself. 

We note for later that the exponential rate of $E(r; n, k, t)$ in these variables 
is given by 
\begin{equation}
\label{eq:exprate}
\limsup_n \frac{\log E(r;n,k,t)}{n} = 
- d_3 \log r + d_1 \log (1 + r) + d_2 \log (1 + 2r).
\end{equation}

Before analyzing the asymptotic behavior of $I(r;n,k,t)$, we discuss the
properties satisfied by the phase term that are essential for the
application of the Laplace approximation. For this and based upon
analytic properties to be clarified shortly, observe first that

\begin{eqnarray}
\label{eq:1st der F at 0}
\frac{\partial F}{\partial\theta}(0;r,d_1,d_2,d_3)&=&i\left\{d_3-\frac{d_1 r}{1+r}-\frac{2d_2 r}{1+2r}\right\},\\
\label{eq:2nd der F at 0}
\frac{\partial^2 F}{\partial\theta^2}(0;r,d_1,d_2,d_3)&=&
\frac{d_1r}{2(1+r)^2}+\frac{d_2r}{(1+2r)^2}.
\end{eqnarray}
Thus, in order for $\theta=0$ to be a stationary point of $F(\theta;r,d_1,d_2,d_3)$, $r$ and $(n,k,t)$ must satisfy the relation $d_3=d_1r/(1+r)+2d_2r/(1+2r)$. A solution $r\ge0$ to this equation is given by the formula
\begin{align}
\label{eq:def r}
r=\frac{2d_3}{(1-3d_3)+\sqrt{(1-3d_3)^2+8d_3(d_1+d_2-d_3)}}.
\end{align}
Note that there is a unique positive solution for $r$.
On the other hand, using that $d_1+2d_2=1$, it follows almost immediately that
\begin{align}
\label{eq:bound 2nd der F}
\frac{\partial^2 F}{\partial\theta^2}(0;r,d_1,d_2,d_3)\ge\frac{r}{2(1+2r)^2}.
\end{align}
Similarly, but after using that $\ln(1-w)\le-w/2$, for all $0\le w\le1$, it follows that
\begin{align}
\label{eq:bound Re F}
\Re\{F(\theta;r,d_1,d_2,d_3)\}\ge\frac{(1-\cos\theta)r}{2(1+2r)^2}.
\end{align}

Thus for given $d_1, d_2, d_3$ in the right range, $F$ has a single
stationary point at $\theta = 0$, satisfying the hypotheses of
Lemma~\ref{thm:Laplace}. Furthermore since $g = 1$ there the Laplace
approximation is uniform as long as no hypotheses change and $r$ remains
bounded away from zero.

In what follows, unless otherwise stated, $r$ always stands for the
short form of the term defined in \eqref{eq:def r}. In addition, we
write $E(n,k,t)$, $I(n,k,t)$ and $F(\theta;d_1,d_2,d_3)$ respectively as
a short form for $E(r;n,k,t)$, $I(r;n,k,t)$ and $F(\theta;r,d_1,d_2,d_3)$.

\subsubsection{The linear case}
\label{sss:linear}

We first study the asymptotic behavior of the coefficient $a(n,k,t)$ for
the regime where $k=\lfloor cn\rfloor$ and $t=\lfloor (n-4)/12\rfloor$,
with $0<c<1/2$. In this case, as $n\to\infty$, $d_1\to5/6$,
$d_2\to1/12$, $d_3\to c$ and $r\to r_c$, where $r_c>0$ is the quantity
defined as 
\begin{align}
r_c:=\frac{2c}{(1-3c)+\sqrt{(1-3c)^2+8c(11/12-c)}}.
\end{align}

In particular, if $c$ is bounded away from zero then for sufficiently large $n$ independent of $c$, $r$ is also bounded away from zero.
Thinking momentarily of $(\theta;r,d_1,d_2,d_3)$ as a vector of unrelated variables, observe that there exists a sufficiently small $0<\delta<\pi$ such that $F(\theta;r,d_1,d_2,d_3)$ is an analytic function of $\theta$ for $|\theta|<2\delta$, for all $(r,d_1,d_2,d_3)$ such that $|r-r_c|<2\delta$.
Thus by Laplace's approximation $I(n,k,t)$ is asymptotically of order
$n^{-1/2}$ as $n \to \infty$. Hence the exponential rate of $p(n, k, t)$
is indeed given by that of $\binom{n}{k}^{-1} E(n, k, t)$. Using
\eqref{eq:binom-approx-linear} and \eqref{eq:exprate} we see that this
rate is 
$$
\frac{5}{6}\ln(1+r_c)+\frac{1}{12}\ln(1+2r_c)-c\ln(r_c)+(1-c)\ln(c)+c\ln(c).
$$ 

It is readily computed that the supremum of the exponential rate for $0
\leq c <11/12$ occurs only when $c\to 0^+$ and has value $0$. Thus
certainly for $0< c < 1/2$, the exponential rate is negative.  This
together with Theorem~\ref{the:1} yields the following result.

\begin{theorem}\label{the:3}
The diameter of a random Cayley digraph of order $n$ and degree $k$ is asymptotically almost surely equal to two provided that $k/n$ remains in a compact subset of the interval $(0,1/2)$ as $n\to\infty$. Furthermore, the convergence is exponentially fast.
\hfill $\Box$
\end{theorem}

Of course the same result will hold for larger values of $c$. Note that
when $c > 11/12$ then for large enough $n$, $k + t > n$ for the value of
$t$ considered here and so $a(n, k, t) = 0$. 

We note in passing that in this case where $k$ is asymptotically linear
in $n$, the analysis of the asymptotics of $a(n, k, t)$ is easily
accomplished by the recently developed methods of Robin Pemantle and
Mark Wilson \cite{p1, p2}. The resulting asymptotic expansion can be
read off almost directly from the explicit expression for $\gf_1$. We
refer to \cite[Section 4.9]{p9} for more details. However, the methods
of \cite{p1, p2} do not work directly in the sublinear case, unlike the
methods of the present paper.

\subsubsection{The sublinear case}
\label{sss:sublinear}

Next we study the asymptotic behavior of $a(n,k,t)$ for the regime where
$k=\lfloor n^\alpha\rfloor$ and $t=\lfloor (n-4)/12\rfloor$, with
$1/2<\alpha<1$. As before, $d_1\to5/6$ and $d_2\to1/12$. However,
$d_3\to 0$ and therefore $r\to 0$ as $n\to\infty$. The new difficulty
here is that the phase term of $I(r;n,k,t)$ converges uniformly to 0 for
all $-\pi\le\theta\le\pi$ as $n\to\infty$. 

To resolve this issue we factor out $r$, 
exploiting the fact that $F(\theta;r,d_1,d_2,d_3)$ is also analytic with respect to
$r$. Indeed, thinking again of $(\theta;r,d_1,d_2,d_3)$ as a vector of
unrelated variables, observe that there exists $\delta>0$ such that
$F(\theta;r,d_1,d_2,d_3)$ is an analytic function of $(\theta;r)$ for
$|\theta|<2\pi$ and $|r|<2\delta$, for all $(d_1,d_2,d_3)$. Thus since
$F(\theta;0,d_1,d_2,d_3)-\frac{\partial F}{\partial\theta}(0;0,d_1,d_2,d_3)\theta=0$, there exists a function
$F_1(\theta;r,d_1,d_2,d_3)$, analytic in $(\theta;r)$ for $|\theta|<2\pi$
and $|r|<2\delta$, such that 
\begin{equation*}
F(\theta;r,d_1,d_2,d_3)-\frac{\partial F}{\partial\theta}(0;r,d_1,d_2,d_3)\theta=r\cdot F_1(\theta;r,d_1,d_2,d_3). 
\end{equation*}

In what follows we write $F_1(\theta;d_1,d_2,d_3)$ as a short form for
$F_1(\theta;r,d_1,d_2,d_3)$. Reverting to our value of $r$ given by \eqref{eq:def r}, we see that 
\begin{equation}
\label{ide:F equal rG}
F(\theta;d_1,d_2,d_3)=r\cdot F_1(\theta;d_1,d_2,d_3).
\end{equation}

Using \eqref{eq:def I(r;n,k,t)} and \eqref{ide:F equal rG}, we obtain

\begin{align}
\label{eq:I with G}
I(n,k,t)=\int_{-\pi}^\pi e^{-nr\cdot F_1(\theta;d_1,d_2,d_3)}d\theta,
\end{align}
for all $(n,k,t)$ such that $0<r\le\delta$. Furthermore, given the
factorization in \eqref{ide:F equal rG}, it follows from \eqref{eq:1st
der F at 0}, \eqref{eq:def r}, \eqref{eq:bound 2nd der F} and
\eqref{eq:bound Re F} that $\frac{\partial
F_1}{\partial\theta}(0;d_1,d_2,d_3)=0$ and that

\begin{align*}
\frac{\partial^2 F_1}{\partial\theta^2}(0;d_1,d_2,d_3)&\ge\frac{1}{2(1+2r)^2},\\
\Re\{F_1(\theta;d_1,d_2,d_3)\}&\ge\frac{1-\cos\theta}{2(1+2r)^2}.
\end{align*}

Now it is readily computed from the definition that 
\begin{equation}
\label{ide:expansion r}
r \sim d_3 - d_3 (d_1 - 5/6) - 2 d_3 (d_2 - 1/12) + 7 d_3^2/6.
\end{equation}

In particular, given that $d_1 \to 5/6$ and $d_2 \to 1/12$, we have 
$n\cdot r\sim n\cdot d_3 = k\to\infty$, as $n\to\infty$.

Since $F_1(\theta;d_1,d_2,d_3)$ is analytic in the disk $|\theta|<2\pi$,
the Laplace approximation can be reapplied but this time to determine
the asymptotic behavior of the integral on the right-hand side of
\eqref{eq:I with G}. It follows that $I(n,k,t)$ is asymptotically of
order $(nr)^{-1/2} \sim k^{-1/2} = O(n^{-1/2})$. As a result, the
exponential growth rate of $p(n,k,t)$ is again given by that of
$\binom{n}{k}^{-1} E(n,k,t)$.

The exponential rate in question is
$$
d_3 \log d_3 + (1 - d_3) \log (1 - d_3) - d_3 \log r + (1 - 2d_2) \log (1 + r) + d_2 \log (1 + 2r).
$$
Using \eqref{ide:expansion r} we see that as $n \to \infty$ this rate is
asymptotic to $-d_3^2/12$. When $k = \lfloor n^{\alpha}\rfloor$ with
$\alpha > 1/2$ the exponential part of $p(n, k, t)$ is therefore
$\exp(n^{(1 - 2\alpha)/12}(1+o(1))$. With the help of Theorem
\ref{the:1} we finally obtain:

\begin{theorem}\label{the:6}
For any constant $\alpha$ such that $1/2<\alpha<1$, the diameter
of a random Cayley digraph of order $n$ and degree $\lfloor
n^{\alpha}\rfloor$ is asymptotically almost surely equal to two.
\hfill $\Box$
\end{theorem}

Note that convergence of the upper bound to zero is faster than 
polynomial, but subexponential.

\subsection{The threshold}
\label{ss:threshold}

We have not yet answered the original question in the introduction,
concerning the threshold for $k = f(n)$ at which the asymptotic value of
$\Pr(\Diam > 2)$ undergoes a phase transition, switching abruptly from $1$
to $0$ as $k$ increases.

The methods above give some useful information on this point. Consider the
simpler analysis of Section~\ref{ss:abelian}, concerning abelian
$2$-groups (similar calculations occur when considering the bounds for
general groups). Assuming that $k = o(t)$, we see that the lower bound
has order of growth equal to that of $b(t, k)$ as $t \to \infty$. As we
have seen, if $k = \Omega(t^\alpha)$ with $\alpha > 1/2$, then $b(t, k)$
converges to zero faster than any polynomial, so the upper bound
converges to zero. To see where the upper bound is asymptotically
constant, we observe from the approximation that $\exp(-\lambda^2 t/4)$
must be of order $t^{-1}$. Thus we require $k \approx 2 \sqrt{ t \log
t}$ for this to occur. At this stage the upper bound converges to $2$;
the more precise $k =  2 \sqrt{t \log t + \log{2}}$ yields a limiting
upper bound of $1$. At this stage the lower bound looks  like $1/t$ and
converges to $0$. The lower bound converges to $1$ only when $k=o(\sqrt{t})$ which 
is not useful. If $k$ grows faster than $\sqrt{t \log t}$,
the upper bound converges to zero. If $k$ grows slower than 
$\sqrt{t\log t}$, then the upper bound goes to infinity with $t$. This gives a
threshold (rather weaker than we might hope for)  in the abelian case. 

In the nonabelian case we can make a similar argument with the upper
bound. However we do not have a good lower bound on the probability. It
may be possible to extract one by refining our arguments of
Sections~\ref{sec:model} and~\ref{sec:estimates}. However Robin Pemantle
(personal communication) has discovered an approach using probabilistic
techniques and along these lines that gives sharper results on the
threshold. Thus we do not proceed further here, preferring to await the
appearance of Pemantle's work.

\section{Conclusions}
\label{sec:conclusions}

We have derived precise information on the event that a random Cayley
digraph has diameter 2, in the abelian group case, and slightly less
precise information in the general case. Many natural questions have
been answered by our asymptotic analysis of upper and lower bounds on
probability. An open question concerns the behaviour in the abelian case
when $k = c \sqrt{n}$. Our upper bound on probability converges to
$\infty$ and the lower bound to $\exp(-c^2/2)$. Perhaps better bounds
will allow us to determine the exact limiting probability using methods
similar to those in this paper.

The genesis of this paper may be of interest. PP and the first JS were
visiting the second JS in Auckland, where they derived the bounds of
Section 3 and posed several questions regarding the asymptotic
behaviour. Their enquiries about asymptotic analysis led from Auckland
to Slovenia (M. Petkovsek) to Pennsylvania (H. Wilf) and then via Robin
Pemantle to ML and MW, the latter being blissfully unaware in Auckland
of the existence of the work going on in the same building!

\bibliographystyle{plain}
\bibliography{cayley}

\end{document}